\documentclass{amsart}

\usepackage{amssymb,amscd,amsthm,amsxtra}
\usepackage{dsfont,latexsym}

\newtheorem{thm}{Theorem}[section]
\newtheorem{corollary}[thm]{Corollary}
\newtheorem{lemma}[thm]{Lemma}
\newtheorem{prop}[thm]{Proposition}
\theoremstyle{definition}
\newtheorem{defn}[thm]{Definition}
\theoremstyle{remark}

\numberwithin{equation}{section}

\newcommand{\eps}{\varepsilon}

\newcommand{\R}{\mathds{R}}
\newcommand{\N}{\mathds{N}}

\newcommand{\CC}{\mathcal{C}}
\newcommand{\FF}{\mathcal{F}}
\newcommand{\EE}{\mathcal{E}}
\newcommand{\KK}{\mathcal{K}}

\usepackage{mathrsfs}
\renewcommand{\mathcal}{\mathscr}

\renewcommand{\le}         {\leqslant}

\renewcommand{\ge}         {\geqslant}

\begin{document}

\author{Ovidiu Savin and Enrico Valdinoci}
\thanks{OS has been supported by
NSF grant 0701037. EV has been
supported
by
FIRB
project ``Analysis and Beyond''
and GNAMPA
project
``Equa\-zio\-ni non\-li\-nea\-ri su va\-rie\-t\`a:
pro\-prie\-t\`a qua\-li\-ta\-tive e clas\-si\-fi\-ca\-zio\-ne
del\-le so\-lu\-zio\-ni''.
Part of this work
was carried out
while EV was visiting Columbia
University.}

\title{$\Gamma$-convergence for
nonlocal phase transitions}

\begin{abstract}
We discuss the $\Gamma$-convergence,
under the appropriate scaling, of the energy functional
$$ \|u\|_{H^s(\Omega)}^2+\int_\Omega W(u)\,dx,$$ with $s \in (0,1)$,
where $\|u\|_{H^s(\Omega)}$ denotes the total contribution from $\Omega$
in the $H^s$ norm of $u$, and $W$ is a double-well potential.

When~$s\in [1/2,\,1)$, we show that the energy
$\Gamma$-converges to the classical minimal surface functional -- while,
when~$s\in(0,\,1/2)$, it is easy to see that the
functional $\Gamma$-converges to the nonlocal minimal surface
functional.
\end{abstract}
\maketitle

\section{Introduction and statement of the main results}

As well-known, the $\Gamma$-convergence, introduced in~\cite{DG,DGF},
is a notion of convergence for
functionals, which tends to be as compatible as possible with the
minimizing features of the energy, and whose limit
is capable to capture essential features of the problem. We refer
to~\cite{DalMaso,Braides} for a
detailed presentation of several basic aspects
and applications of $\Gamma$-convergence; see also~\cite{tartar}
for applications to homogenization theory.

Making it possible to study
the asymptotics of variational problems indexed by a parameter, the
$\Gamma$-convergence has become a standard tool in dealing
with singularly perturbed energies as the ones
arising in the theory of phase transitions (see~\cite{MM}),
where the dislocation energy of a double well potential~$W$ is compensated
by a small gradient term which avoids the formation of unnecessary
interfaces, leading to a total energy which is usually written as
\begin{equation}\label{old}
\int {\eps^2} |\nabla u|^2+W(u) \, dx,\end{equation}
with~$\eps\rightarrow 0^+$.

The purpose of this paper is to develop a $\Gamma$-convergence
theory for a nonlocal analogue of the energy above, in which
the gradient term in~\eqref{old}
is replaced by a fractional, Gagliardo-type, (semi)norm
of the form~$\eps^{2s} \|u\|_{H^s}^2$, with~$s \in (0,1)$
(see below for precise definitions and statements).
Notice that, formally, the gradient term in~\eqref{old}
corresponds to the case~$s=1$.

The study of such a nonlocal contribution is quite
important for the applications, since the classical gradient
term takes into account the interactions at small scales between
the particles of the medium, but loses completely the long scale
interactions. In this spirit, it is relevant to know
whether or not the $\Gamma$-limit of the functional
is local -- that is, whether or not the long range interactions
affect the limit interface.

{F}rom the point of view of the pure mathematics, nonlocal problems
are also relevant because new techniques are usually needed
to understand and estimate the contributions coming from far.
We refer, in particular, to~\cite{CRS} for the definition and
the basic features of nonlocal minimal surfaces, which are
the natural analogue of the classical sets of minimal
perimeter (as in~\cite{Giusti}).
In fact, we will show that the $\Gamma$-limit of
our functional will be the standard minimal surface functional
when~$s\in [1/2,1)$ and the nonlocal one when~$s\in (0,1/2)$.

Now, we introduce the formal setting in which we work.
We consider a bounded domain~$\Omega\subset\R^n$, $n\ge2$,
with 
complement~$\CC\Omega$.
We define
$$ X:=\big\{ u\in L^\infty(\R^n) \,:\; \| u\|_{L^\infty(\R^n)}\le 1
\big\},$$
the space of admissible functions $u$.
We say that a sequence~$u_k\in X$ converges to~$u$ in~$X$
if~$u_k$ converges to~$u$ in~$L^1_{\rm loc}(\R^n)$.

We define
$$ \KK(u,\Omega):=\frac{1}{2}\int_\Omega\int_\Omega\frac{
|u(x)-u(y)|^2
}{
|x-y|^{n+2s}}\,dx\,dy
+\int_\Omega\int_{\CC\Omega}\frac{
|u(x)-u(y)|^2
}{
|x-y|^{n+2s}}\,dx\,dy,$$
the $\Omega$ contribution in the $H^s$ norm of $u$ $$\int_{\R^n}
\int_{\R^n} \frac{
|u(x)-u(y)|^2
}{
|x-y|^{n+2s}}\,dx\,dy,$$ i.e. we omit the set where $(x,y) \in \CC\Omega
\times \CC \Omega$ since all $u\in X$ are fixed outside $\Omega$.

The energy functional $J_\eps$ in $\Omega$ is defined as
$$J_\eps(u,\Omega):=\eps^{2s}\KK(u,\Omega) + \int_\Omega W(u) \, dx.$$
Such functional may be seen as the nonlocal analogue
of the classical one in~\eqref{old}.

Throughout the paper we assume that $W:[-1,1] \to [0, \infty)$,
\begin{equation}\label{Wcond}\begin{split}
& W \in C^2([-1,1]), \quad W(\pm 1)=0, \quad W>0 \quad \mbox{in
$(-1,1)$},\\
&\qquad
W'(\pm 1)=0 \quad{\mbox{ and }}\quad  W''(\pm
1)>0.\end{split}\end{equation}

We remark that, differently from several nonlocal models
considered in the literature (see e.g. \cite{ABCPJSP,ABS,GM2}
and references therein),
we deal with an arbitrarily large number of space dimensions,
no periodicity in space is assumed, and we consider the 
full interaction
among all the space $\Omega$ versus $\R^n$
(i.e., from the physical point of
view, the
particles in the domain $\Omega$ interact with the ones in the whole of
the space
$\R^n$, not only with the ones in $\Omega$).

Since~$\Gamma$-convergence is expecially designed
for minimizers, we recall the following notation:

\begin{defn}\label{D}
We say that $u$ is a minimizer for $J_\eps$ in 
an open, possibly unbounded, set~$\Omega\subset\R^n$
if, for any open subset~$U$
compactly included in $\Omega$, we have that 
$$J_\eps(u,U) < \infty,$$ and
$$J_\eps(u,U) \le J_\eps(v,U)$$ for any $v$ which coincides with
$u$ in $\CC U$.
\end{defn}

It is worth to notice that
if $u$ minimizes $J_\eps$ in
$\Omega$ then it minimizes
$J_\eps$ in any subdomain $\Omega'\subset \Omega$.

We deal with the functional~$\FF_\eps:X\rightarrow\R\cup\{+\infty\}$
defined as
\begin{equation*}\begin{split}&
\FF_\eps(u,\Omega) :=\left\{
\begin{matrix}
\eps^{-2s} J_\eps(u, \Omega) & {\mbox{ if $s\in(0,\,1/2)$,}} \\
\ \\
|\eps \log \eps|^{-1} J_\eps(u, \Omega) & {\mbox{ if $s=1/2$,}}\\
\ \\
\eps^{-1} J_\eps(u, \Omega) & {\mbox{ if $s\in(1/2,\,1)$.}}
\end{matrix}
\right.\end{split}\end{equation*}
The functional~$\FF_\eps$ may be seen as the ``right''
scaling of~$J_\eps$, i.e. the one
that comes from the dialation invariance of
the space and that possesses 
an interesting~$\Gamma$-limit, in relation with phase transitions.

In the case when $s\in (0,1/2)$, the limiting  functional~$\FF:X
\rightarrow\R\cup\{+\infty\}$
is defined as
\begin{equation}\label{FF1}
\FF(u,\Omega):=\left\{\begin{matrix}
\KK(u,\Omega)
& {\mbox{ if $u|_\Omega = \chi_E -\chi_{\CC E}$, for some set $E\subset
\Omega$}}
\\+\infty & {\mbox{otherwise.}}
\end{matrix}\right.
\end{equation}
In this case,~$\FF$
agrees with the nonlocal area functional of $\partial E$ in $\Omega$ that
was studied in~\cite{CRS, CV, ADP}.
Remarkably, such nonlocal area functional is well defined
exactly when~$s\in(0,\,1/2)$.

In the case when $s \in [1/2,1)$ the limiting functional~$\FF:X\rightarrow\R\cup\{+\infty\}$
is defined as
\begin{equation}\label{FF2}
\FF(u, \Omega):=\left\{\begin{matrix}
c_\star {\rm Per} (E,\Omega)
& {\mbox{ if $u|_\Omega = \chi_E -\chi_{\CC E}$, for some set $E\subset
\Omega$}}
\\+\infty & {\mbox{otherwise,}}
\end{matrix}\right.
\end{equation}
where $c_\star$ is a constant depending on $n$, $s$ and $W$, which
will be explicitly determined in the sequel, in dependence of
a suitable $1$D minimal profile (see Theorem \ref{t1}
and \eqref{DEFCSTAR} for details).

Here above and in the rest of the paper, we use the standard
notation ${\rm Per}(E,U)$
to denote the perimeter of a set $E$ in an open set $U\subseteq
\R^n$ (see, e.g., \cite{Giusti}).

Then, the results we prove here are the following:

\begin{thm}\label{TH1}
Let~$s\in(0,\,1)$. Then,
$\FF_\eps$ $\Gamma$-converges to~$\FF$, i.e.,
for any~$u\in X$,
\begin{itemize}
\item[{(i)}] for any~$u_\eps$ converging to~$u$ in~$X$,
$$ \FF(u,\Omega)\le\liminf_{\eps\rightarrow0^+}\FF_\eps(u_\eps,\Omega),$$
\item[{(ii)}] if $\Omega$ is a Lipschitz domain,
there exists~$u_\eps$ converging to~$u$
in~$X$ such that
$$ \FF(u,\Omega)\ge\limsup_{\eps\rightarrow0^+}\FF_\eps(u_\eps,\Omega).$$
\end{itemize}
\end{thm}

\begin{thm}\label{comp}
If $\FF_\eps(u_\eps, \Omega)$ is uniformly bounded for a sequence of $\eps
\to 0^+$, then there exists a convergent subsequence 
\begin{equation}\label{conv}u_\eps \rightarrow
u_*:=\chi_E-\chi_{\CC E} \quad \mbox{ in $L^1(\Omega)$.}
\end{equation} Moreover,
let~$u_\eps$ minimize $\FF_\eps$ in $\Omega$:
\begin{itemize}
\item[{(i)}] if~$s\in(0,1/2)$ and $u_\eps$ converges weakly to $u_o$
in $\CC
\Omega$, then $u_*$ minimizes $\FF$ in~\eqref{FF1} among all the 
functions that coincide with
$u_o$ in $\CC\Omega$;
\item[{(ii)}] if~$s\in[1/2,1)$, then $u_*$ minimizes $\FF$ in
\eqref{FF2}.
Also, for any open set~$U \subset \subset
\Omega$ we have $$\limsup_{\eps\rightarrow0^+} \FF_\eps(u_\eps,U) \le c_*
{\rm Per}\,(E,\overline{U}).$$
\end{itemize}
\end{thm}

As we will see in the rest of the paper, the $\Gamma$-convergence
for~$s\in(0,1/2)$ is elementary, so
the hard task is to deal with the case~$s\in[1/2,1)$.
Naifly, this reflects the fact that in such a case
the nonlocal charachter of the problem
gets localized in the $\Gamma$-limit, hence the estimates
need to carefully take into account the nonlocal contributions
and their local counterparts, and balance the ones with the 
others in a precise way.

We recall that there are several results available
in the literature concerning the approximation
of the perimeter with nonlocal functionals.
As far as we understand, all these results are
related to our
Theorems \ref{TH1}
and \ref{comp}
(as well as to each other),
but their statements are quite different from ours and the proofs
are based on different techniques.
In particular, we recall \cite{ABSCRAS}, which considered
a $H^{1/2}$ norm inside a
one-dimensional domain with no contribution
coming from the outside. As remarked to us by \cite{ABSper},
the extension of the results
in \cite{ABSCRAS} to higher dimension is implicitly contained
in \cite{ABS}, though not explicitly mentioned.
Moreover,
in \cite{GM1,GM2} the $\Gamma$-convergence
of a functional driven by a norm of type $H^{1/2}$
and a more complicated potential on a
two-dimensional
square or torus, under a suitable pinning condition,
was studied in detail.

Also, in \cite{ABCPJSP,ABEJAM}, the $\Gamma$-convergence of an interaction
energy with a double integral weighted by a summable kernel is
considered: here, we take into account integrands with
a more severe singularity so that many technical difficulties
need to be overcome.

{F}rom the results in
Theorem \ref{TH1} and
\ref{comp}, 
it is also 
possible
to have optimal estimates on the width of the asymptotic
interface of minimizers.
Indeed, in~\cite{density paper} we proved the following energy bound and uniform density estimate for minimizers of $\FF_\eps$.

\begin{thm}\label{EB} If $u_\eps$ minimizes $\FF_\eps$ in $B_{1+2\eps}$ then
$$\FF_\eps(u_\eps,B_1) \le \overline{C},$$
with $\overline{C}$ depending on $n$, $s$, $W$.
\end{thm}

\begin{thm}\label{DE}
If $u_\eps$ minimizes $\FF_\eps$ in $B_r$ and $u_\eps(0)>\theta_1$ then
$$| \{ u_\eps>\theta_2\}\cap B_r|\,\ge\, \overline{c} \,r^n$$ provided that $\eps \le c(\theta_1,\theta_2) r$,
where $\overline c>0$ depends only on $n$, $s$, $W$ 
and~$c(\theta_1,\theta_2)>0$ depends also on $\theta_1, 
\theta_2\in(-1,1)$. \end{thm}

As a consequence of these theorems we obtained in~\cite{density paper} that the convergence in \eqref{conv} is better when dealing with minimizers. 
More precisely, we showed that the level sets of minimizers $u_\eps$ of $\FF_\eps$ converge locally uniformly to $\partial E$.

For the proof of Theorems~\ref{EB} and~\ref{DE}, see~\cite{density paper}.
We also refer to~\cite{ABS,GP,Gon}, where other types of nonlocal
models have been considered (in particular, a three-dimensional fluid
with boundary and weight inhomogeneity
of distance type, whose energy bounds the Gagliardo norm, see
Theorem~19 in~\cite{Gon}).

The proof of Theorems~\ref{TH1}
and~\ref{comp} when~$s\in(0,\,1/2)$ is elementary and it
is contained in Section~\ref{SS1}.
In Section \ref{4}
we prove the compactness needed in Theorem \ref{comp} in the case $s \ge
1/2$.
In Section \ref{SS2} we
prove Theorem~\ref{TH1}
and Theorem ~\ref{comp} (ii) when~$s\in [1/2,1)$ by
interpolating the functions candidate to the minimization.
For this, a careful analysis on the energy contribution
across the gluing of the interpolation is needed, as well as
some measure theoretic result of~\cite{density paper}. 

Several arguments in the sequel will be based on some
preliminary considerations, whose detailed proofs can be found
in~\cite{PSV}.

Finally, we conclude the introduction with a notation that will be used throughout the paper. For simplicity we denote

\begin{equation}\label{p.6}u(E,F):=\int_E\int_F 
\frac{|u(x)-u(y)|^2}{|x-y|^{n+2s}}\, dx dy.
\end{equation}
Clearly,
$u(E,F)=u(F,E)$, and if $E_1$ and $E_2$ are disjoint, then $$u(E_1\cup E_2,F)=u(E_1,F)+u(E_2,F).$$
Using this notation, the $\Omega$ contribution in the $H^s$ norm of $u$ can be written as
$$\KK(u, \Omega)=\frac12 u(\Omega, \Omega)+u(\Omega, \CC \Omega).$$

\section{Proof of Theorems~\ref{TH1}
and~\ref{comp} when~$s\in(0,\,1/2)$}\label{SS1}

Throughout this section we assume~$s\in(0,\,1/2)$.

{\it Proof of Theorem~\ref{TH1}.} Recalling~\eqref{FF1},
we observe that
\begin{equation}\label{FK}
{\mbox{if $u\Big|_\Omega=\chi_E-\chi_{\CC E}$, then
$\FF_\eps(u,\Omega)=\FF(u,\Omega)=\KK(u,\Omega)$.}}
\end{equation}
Now, we
prove~(i).
For this, let~$u_\eps$ converging to~$u$ in~$X$.
If
$$ \liminf_{\eps\rightarrow0^+}\FF_\eps(u_\eps,\Omega)=+\infty,$$
then~(i) is obvious, so we may suppose that
$$ \liminf_{\eps\rightarrow0^+}\FF_\eps(u_\eps,\Omega)=\ell<+\infty.$$
We take a subsequence, say~$u_{\eps_k}$ attaining the above
limit.

Then, we take a further subsequence, say~$u_{\eps_{k_j}}$,
that converges to~$u$ almost everywhere.
Therefore,
$$ \ell=\lim_{k\rightarrow+\infty}\FF_{\eps_k}(u_{\eps_k},\Omega)
=\lim_{j\rightarrow+\infty}\FF_{\eps_{k_j}}(u_{\eps_{k_j}},\Omega)
\ge \lim_{j\rightarrow+\infty}
\frac{1}{\eps^{2s}_{k_j}} \int_\Omega W(u_{\eps_{k_j}}(x))\,dx.$$
Consequently,
$$ \int_\Omega W(u(x))\,dx =\lim_{j\rightarrow+\infty}
\int_{\Omega} W(u_{\eps_{k_j}}(x))\,dx = 0.$$
This implies that~$u(x)\in\{-1,+1\}$ for almost any~$x\in\Omega$,
that is,~$u\Big|_\Omega=\chi_E-\chi_{\CC E}$ for a suitable set~$E$. And
so, by Fatou Lemma
and~\eqref{FK}, we conclude that
$$ \liminf_{\eps\rightarrow0^+}\FF_\eps(u_\eps,\Omega)\ge
\liminf_{\eps\rightarrow0^+}\KK (u_\eps,\Omega)\ge \KK(u,\Omega)=\FF(u,\Omega),$$
proving~(i). Now, we prove~(ii).

For this, we may suppose that~$u\Big|_\Omega=\chi_E-\chi_{\CC E}$ for a
suitable set~$E$, otherwise~(ii)
is obvious. Then, we choose~$u_\eps:=u$ and we use~\eqref{FK}
to see that~$\FF_\eps(u_\eps,\Omega)=\FF(u,\Omega)$,
which obviously implies~(ii). This completes the proof of
Theorem~\ref{TH1}.

\

{\it Proof of Theorem \ref{comp}.} Since $s\in(0,\,1/2)$, the uniform 
bound on $\FF_\eps$ gives a uniform bound of the Gagliardo norm 
$\KK(u_\eps, \Omega)$,
and the compactness claim in~\eqref{conv}
is quite standard, see for 
example Section~6 in~\cite{PSV}. It remains to prove (i).

As a result of
Definition~\ref{D}, it suffices to consider the case 
when $\Omega$ is bounded and smooth.
In this case, one has that
\begin{equation}\label{SXC}
\int_{\CC
\Omega}\int_\Omega
\frac{2}{|x-y|^{n+2s}}dxdy<\infty
.\end{equation}

Let~$v\in X$ be an arbitrary function with$$v\Big|_\Omega=\chi_F-\chi_{\CC F}$$
for some set~$F$, and ~$v=u_o$
in~$\CC\Omega$. For any~$y\in\CC\Omega$, let
$$ \psi(y):=\int_\Omega\frac{v(x)}{|x-y|^{n+2s}}\,dx\quad{\mbox{
and
}} \quad\Psi(y):=\int_\Omega\frac{u_*(x)}{|x-y|^{n+2s}}\,dx\,,$$
where~$u_*$ is as in~\eqref{conv}.
We remark that $\psi(y)$ and $\Psi(y)$ are in $L^1(\CC \Omega)$, since
$$\int_{\CC \Omega}|\psi(y)|+|\Psi(y)|
\, dy \le \int_{\CC 
\Omega}\int_\Omega 
\frac{2}{|x-y|^{n+2s}}dxdy 
< \infty,$$
thanks to~\eqref{SXC}.

By the weak convergence of~$u_\eps$, and the fact that $|u_\eps|$
and~$|u_o|$ 
are uniformly bounded
\begin{equation}\label{weak1}
\lim_{\eps\rightarrow 0^+} \int_{\CC\Omega}
\big( u_\eps(y)-u_o(y)\big)\phi(y)\,dy\,=\,0\,,\end{equation} for any 
$\phi \in L^1(\R^n)$.
Moreover, by the strong convergence of $u_\eps$ in $\Omega$,
\eqref{SXC},
and the Dominated Convergence Theorem, we have that
\begin{equation}\label{d67f8t1515299}\begin{split}
\lim_{\eps\rightarrow0^+} \int_{\Omega} |u_\eps(x)|^2
\int_{\CC\Omega}\frac{
dy}{|x-y|^{n+2s}}\,dx
=
\int_{\Omega} |u_*(x)|^2
\int_{\CC\Omega}\frac{
dy}{|x-y|^{n+2s}}\,dx
\end{split}\end{equation}
and
\begin{equation}\label{d67f8t1515299bis}\begin{split}
& \lim_{\eps\rightarrow0^+} \left| \int_{\Omega} \big(
u_\eps(x)-u_*(x)\big) \int_{\CC\Omega}\frac{
u_\eps(y)\,dy}{|x-y|^{n+2s}}\,dx \right|
\\ &\qquad\quad\le\,
\lim_{\eps\rightarrow0^+} \int_{\Omega} \big|
u_\eps(x)-u_*(x)\big| \int_{\CC\Omega}\frac{
dy}{|x-y|^{n+2s}}\,dx
\,=\,0.\end{split}
\end{equation}
On the other hand,
making use of the notation in~\eqref{p.6}, we deduce from
Fatou Lemma that
\begin{equation}\label{d67f8t1515299tris}
\liminf_{\eps\rightarrow0^+} u_\eps(\Omega,\Omega)\ge 
u_*(\Omega,\Omega).
\end{equation}
Let also
$$ v_\eps(x):=\left\{\begin{matrix}
v(x) & {\mbox{ if $x\in\Omega$,}}\\
u_\eps (x) & {\mbox{ if $x\in\CC\Omega$.}}
\end{matrix}\right.$$Recalling that~$u_\eps$
is minimal, we obtain
that
\begin{eqnarray*}
0 &\le& \FF_\eps(v_\eps,\Omega)-\FF_\eps(u_\eps,\Omega)
\\ &=&\KK(v_\eps, \Omega)-\KK(u_\eps,\Omega)-\eps^{-2s}
\int_\Omega W(u_\eps(x))\,dx
\\ &\le&\frac{1}{2}\Big( v(\Omega,\Omega)-u_\eps(\Omega,\Omega)\Big)
+\int_\Omega\left(\int_{\CC\Omega}\frac{
|v(x)-u_\eps(y)|^2-|u_\eps(x)-u_\eps(y)|^2
}{|x-y|^{n+2s}}\,dy\right)\,dx
\\ &=&
\frac{1}{2}\Big( v(\Omega,\Omega)-u_\eps(\Omega,\Omega)\Big)
\\&&\quad+\int_\Omega |v(x)|^2
\int_{\CC\Omega}\frac{dy}{|x-y|^{n+2s}}\,dx-
\int_\Omega |u_\eps(x)|^2 \int_{\CC\Omega}\frac{dy}{|x-y|^{n+2s}}
\\ &&\quad+2\int_{\CC\Omega} u_\eps(y)\Psi(y)\,dy-
2\int_{\CC\Omega} u_\eps(y)\psi(y)\,dy
\\ &&\quad+2
\int_{\Omega} \big(
u_\eps(x)-u_*(x)\big) \int_{\CC\Omega}\frac{
u_\eps(y)\,dy}{|x-y|^{n+2s}}\, dx.
\end{eqnarray*}
Consequently, recalling that~$v(y)=u_*(y)=u_o(y)$ for any~$y\in\CC\Omega$
and using ~\eqref{weak1},
\eqref{d67f8t1515299},
\eqref{d67f8t1515299bis}
and~\eqref{d67f8t1515299tris}, we obtain
\begin{eqnarray*}
0&\le&
\frac{1}{2}\Big( v(\Omega,\Omega)-u_*(\Omega,\Omega)\Big)
\\ &&\quad+\int_\Omega |v(x)|^2
\int_{\CC\Omega}\frac{dy}{|x-y|^{n+2s}}\,dx-
\int_\Omega |u_*(x)|^2 \int_{\CC\Omega}\frac{dy}{|x-y|^{n+2s}}
\\&&\quad+2\int_{\CC\Omega} u_o (y)\Psi(y)\,dy-
2\int_{\CC\Omega} u_o(y)\psi(y)\,dy\\ &=&
\frac{1}{2}\Big( v(\Omega,\Omega)-u_*(\Omega,\Omega)\Big)
+\int_\Omega\left(\int_{\CC\Omega}\frac{
|v(x)-v(y)|^2-|u_*(x)-u_*(y)|^2
}{|x-y|^{n+2s}}\,dy\right)\,dx
\\ &=& \FF(v,\Omega)-\FF(u_*,\Omega)
.\end{eqnarray*}
This proves
claim~(i) of
Theorem~\ref{comp}, and it ends the proof of
Theorem~\ref{comp}.

\section{Compactness for $s \ge 1/2$}\label{4}

Here, we prove the compactness claimed in
Theorem \ref{comp}  when $s\in [1/2,\,1)$
(and this range of~$s$
will be assumed throughout this section).
An important tool for our estimate is Proposition~4.3 of~\cite{density 
paper}, which provides a lower bound for the double integral
$$L(A,D):=\int_A \int_D \frac{1}{|x-y|^{n+2s}}dxdy.$$
For the convenience of the reader we state it below.
\begin{prop}\label{P1}
Let $ s\in [1/2,1)$. Let $A$, $D$ be disjoint subsets of a cube $Q 
\subset \R^n$ with 
\begin{equation}\label{bad0}
\min\{|A|,|D|\} \ge \sigma |Q|,
\end{equation}
for some $\sigma 
>0$. Let $B=Q\setminus(A\cup D)$. Then,
$$ L(A,D)\ge \left\{
\begin{matrix}
\delta |Q|^\frac{n-1}{n}\,\log (|Q|/{|B|}) &\qquad{\mbox{ if $s=1/2$,}}\\
\ \\
\delta |Q|^\frac{n-2s}{n}\,(|Q|/|B|)^{2s-1}&\qquad{\mbox{ if 
$s\in(1/2,1)$.}}\\
\end{matrix}
\right.$$
with $\delta >0$ depending on $\sigma$, $n$ and $s$.
\end{prop}

Also, it is
convenient to define \begin{equation}\label{IE}
I_\eps(u, \Omega)=\left\{
\begin{matrix}
\displaystyle\frac{1}{ 2|\log \eps|} u(\Omega,\Omega)+
\displaystyle\frac{1}{\eps|\log
\eps|}\int_\Omega W(u)\,dx & {\mbox{ if $s=1/2$,}}\\
\ \\
\displaystyle\frac{\eps^{2s-1}}{2} u(\Omega, \Omega) +
\displaystyle\frac{1}{\eps}\int_\Omega W(u)\,dx & {\mbox{ if
$s\in(1/2,\,1)$.}}
\end{matrix}
\right.\end{equation}
Notice that~$I_\eps(u_\eps,\Omega)$ depends only on the values of $u$ in
$\Omega$. We list some useful properties of $\FF_\eps$ and $I_\eps$
that follow immediately from their definition:

a) $I_\eps$ is bounded by $\FF_\eps$, i.e.
$$\FF_\eps(u,\Omega) \ge
I_\eps(u, \Omega),$$

b) $\FF_\eps$ is subadditive, i.e. if $E$ and~$F$ are disjoint sets then
$$\FF_\eps(u,E \cup F) \le \FF_\eps(u,E) + \FF_\eps(u,F),$$

c) $I_\eps$ is superadditive, i.e. if $E$ and~$F$ are disjoint sets
then$$I_\eps(u, E\cup F) \ge I_\eps(u,E)+I_\eps(u,F). $$
\medskip

As a consequence of~\eqref{IE} and Proposition~\ref{P1},
we obtain
\begin{lemma}\label{help}
Let~$\sigma\in(0,1/4)$ and~$u\in X$. Let $Q$ be a cube in $\R^n$. If
\begin{equation}\label{bad}
|\{u\ge 1-\sigma\} \cap Q| \ge \sigma |Q|
\;{\mbox{ and }}\;
\quad |\{u\le -1+\sigma\} \cap 
Q| \ge \sigma |Q|
\end{equation}
then, for all small $\eps$
\begin{equation}\label{badg}
I_\eps (u,Q) \ge c(\sigma)|Q|^\frac{n-1}{n}.
\end{equation}
where $c(\sigma)>0$ depends on $\sigma$ and on $n$, $s$, $W$.
\end{lemma}
\begin{proof}
Define $$A:=\{u\ge 1-\sigma\} \cap Q, \quad D:=\{u\le
-1+\sigma\} \cap Q,$$ 
$$ B:=\{|u|\le 1-\sigma\} \cap Q=Q\setminus(A\cup D).$$
If $$|B| \ge \left \{ \begin{matrix}
\eps |\log \eps||Q|^\frac{n-1}{n}, \quad \mbox{and $s=1/2$}\\
\ \\
\eps |Q|^\frac{n-1}{n}, \quad \mbox{and $s>1/2$},
\end{matrix}
\right.
$$
then the potential energy in $I_\eps(u,Q)$ satisfies \eqref{badg} for some 
small $c(\sigma)$, and there is nothing to prove.
Otherwise we apply Proposition~\ref{P1},
noticing that~\eqref{bad0} is satisfied because of~\eqref{bad}:
we obtain 
$$u(Q,Q) \ge u(A,D) 
\ge 
L(A,D) \ge \left \{ \begin{matrix}
\delta(\sigma) \log \left(\displaystyle\frac{|Q|^\frac 1n}{\eps \log \eps}
\right) |Q|^\frac{n-1}{n}, \quad \mbox{and $s=1/2$}\\
\ \\
\delta(\sigma)\eps^{1-2s} |Q|^\frac{n-1}{n}, \quad \mbox{and $s>1/2$.}
\end{matrix}
\right.
$$
This shows that the kinetic energy in $I_\eps(u,Q)$ satisfies \eqref{badg} 
provided that, in the case $s=1/2$, $\eps \le \eps_0(|Q|)$.
\end{proof}
Here is the compactness needed for Theorem \ref{comp}:
\begin{prop}\label{7477.4}
Let~$\Omega$ be an open, bounded subset of~$\R^n$
and~$u_\eps \in X$, with~$\eps>0$.

If $$\liminf_{\eps \to 0^+} \FF_\eps (u_\eps,\Omega)<+\infty,$$
then $u_\eps$ has a subsequence converging in~$L^1(\Omega)$ to $\chi_E 
-\chi_{\CC E}$, for a suitable~$E\subseteq\R^n$. Moreover, $${\rm Per}(E, 
\Omega) < \infty.$$
\end{prop}
\begin{proof}
We prove that the
set~$u_\eps$ is
totally
bounded in
$L^1(\Omega)$,
i.e. for any $\delta>0$
there exists a finite set~${\mathcal{S}}\subset L^1(\Omega)$
such that for any small $\eps$ there exists~$\psi_\eps \in {\mathcal{S}}$ 
with
\begin{equation}\label{asdiff}
\| u_\eps-\psi_\eps\|_{L^1(\Omega)} \le \delta.
\end{equation}
By passing if necessary to a subsequence we assume
\begin{equation}\label{C_0}
C_0 \ge \FF_\eps (u_\eps,\Omega),
\end{equation} for some constant $C_0$. Fix $\sigma>0$ small.
We decompose the space in cubes $Q_i$ of size $\rho$ with $\rho>0$ small, 
depending on $\sigma$ and $\delta$, to be made precise later. Let 
$$K:=\bigcup_{Q_i \subset \Omega} Q_i $$
denote the collection of these cubes which are included in $\Omega$. We 
decompose $K$ in three sets $K_+$, $K_-$, $K_0$ as follows
\begin{align*}
K_+:=& \bigcup_{Q_i \subset F_+}Q_i, \quad \quad F_+=
\Big\{Q_i \in K
{\mbox{ s.t. }} 
|\{u_\eps<-1+\sigma\}\cap Q_i| < \sigma |Q_i|\Big\},\\
K_-:=&\bigcup_{Q_i \in F_-} Q_i, \quad \quad F_-:=
\Big\{Q_i \in K\setminus 
K_+{\mbox{
s.t. }}|\{u_\eps>1-\sigma\}\cap Q_i| < \sigma |Q_i|\Big\},\\
K_0:=& K \setminus (K_+\cup K_-).
\end{align*}
We define $\psi_\eps$ to be 1 in $K_+$, and $-1$ otherwise.
If $\rho$ is sufficiently small then 
\begin{equation} \label{XV1} |\Omega \setminus K| \le \delta 
/8.  \end{equation}
We have $$C_0 \ge \FF_\eps(u_\eps,\Omega) \ge I_\eps(u_\eps,K_0) \ge 
\sum_{Q_i \subset K_0} I_\eps(u_\eps,Q_i) \ge 
\frac{|K_0|}{\rho^n}c(\sigma) \rho^{n-1},$$
where in the last inequality we used Lemma \ref{help}. Hence
\begin{equation} \label{XV2}|K_0| \le 
C(\sigma,C_0) \rho \le  \delta / 8,\end{equation}
provided that $\rho$ is 
small 
enough.

{F}rom \eqref{C_0} we also see that for all small $\eps$
$$|\{|u_\eps | \le 
1-\sigma\} 
\cap \Omega| \le C(\sigma) \int_\Omega W(u_\eps) \, dx \le C(\sigma,C_0) 
\eps^{1/2} \le \delta / 8.$$
Therefore
\begin{equation}\label{XV3}
\int_{
\{|u| \le
1-\sigma\}
\cap \Omega
}|u_\eps-\psi_\eps |\, dx\le 2 |\{|u| \le
1-\sigma\}
\cap \Omega|\le\delta/4.
\end{equation}
Moreover,\begin{eqnarray*}
&& |K_+\cap \{ u_\eps<-1+\sigma\}|=\sum_{Q_i\subset K_+}
|Q_i\cap \{ u_\eps <-1+\sigma\}|\\
&&\qquad<\sigma\sum_{Q_i\subset K_+}
|Q_i|=\sigma |K_+|
\end{eqnarray*}
and so
\begin{equation}\label{XX10}
\int_{K_+\cap 
\{u_\eps<-1+\sigma\}
}|u_\eps-\psi_\eps |\, dx\le
2 |K_+\cap \{ u_\eps<-1+\sigma\}|\le 2\sigma |K_+|
.\end{equation}
In the same way, we obtain
\begin{equation}\label{XX11}
\int_{K_- \cap
\{u_\eps >1-\sigma\}
}|u_\eps-\psi_\eps |\, dx\le
2\sigma |K_-|.\end{equation}
On the other hand,
\begin{equation}\label{ER0}
\begin{split}
& \int_{K_-\cap
\{u_\eps <-1+\sigma\}
}|u_\eps-\psi_\eps |\, dx+
\int_{K_+\cap
\{u_\eps >1-\sigma\}
}|u_\eps-\psi_\eps |\, dx\\
=\,& \int_{K_-\cap  
\{u_\eps <-1+\sigma\}
}|u_\eps+1|\, dx+
\int_{K_+\cap 
\{u_\eps >1-\sigma\}
}|u_\eps-1 |\, dx\\
\le\,& \sigma |{K_-\cap
\{u_\eps <-1+\sigma\}
}|+\sigma |{K_+\cap
\{u_\eps >1-\sigma\}
}|
\\ \le\,& \sigma |K_+\cup K_-|.
\end{split}
\end{equation}
{F}rom~\eqref{XX10}, \eqref{XX11}
and~\eqref{ER0}, we conclude that
\begin{equation*}
\int_{(K_-\cup K_+)\cap
\{|u_\eps| >1-\sigma\}
}|u_\eps-\psi_\eps |\, dx\le 3\sigma|K_+\cup K_-|.
\end{equation*}
This and~\eqref{XV3} yield that
$$ \int_{K_-\cup K_+   
}|u_\eps-\psi_\eps |\, dx
\le\delta/2$$
as long as~$\sigma$ is small enough.

{F}rom the latter inequality and the ones in~\eqref{XV1},
and~\eqref{XV2} we obtain
\begin{align*}
\int_{\Omega}|u_\eps-\psi_\eps |\, dx &\le 2|\Omega \setminus K| 
+2|K_0|+\int_{K_+\cup K_-} |u_\eps-\psi_\eps |\, dx \le \delta.
\end{align*}
The set $S$ of all $\psi_\eps$ is clearly finite and our claim is proved. 
Since $|\psi_\eps|\equiv 1$ we can easily conclude that there exists a 
convergent subsequence of $u_\eps$'s
 in $L^1(\Omega)$ to a function of the form $\chi_E-\chi_{\CC E}$ for some
set $E$.
It remains to show that if $u_\eps$ converges to $\chi_E-\chi_{\CC E}$ 
then $E$ has finite perimeter in $\Omega$. As above, we decompose $\R^n$ 
into cubes $Q_i$ of size $\rho$ and define 
$$\phi_\rho=\left\{\begin{array}{l}
 1\quad \quad \mbox{in $Q_i$ if $|E\cap Q_i| \ge 1/2 |Q_i|$} \\
 -1 \quad \quad \mbox{otherwise.}
 \end{array}
 \right.
$$ We also define $$\tilde \phi_\rho:=\phi_\rho * g_\rho$$ where $g_\rho$ 
is a mollifier defined in $B_\rho$, and we remark that $$|\nabla \tilde 
\phi_\rho| \le C/\rho.$$
{F}rom Lebesgue Theorem, $\psi_\rho$ and $\tilde \psi_\rho$ converge to 
$\chi_E-\chi_{\CC E}$ as $\rho \to 0^+$.
Now we estimate the $BV$ norm of $\tilde \psi_\rho$ by counting the number 
of cubes
$Q_i$ in $\Omega$ at distance greater than $\sqrt n \rho$ from $\partial 
\Omega$, i.e. $Q_i \in \Omega_{\sqrt n \rho}$,
for which $\tilde \psi_\rho$ is not constant (1 or $-1$) in $Q_i$. Denote 
the set of such cubes by $F$.
If $Q_i \in F$, then the cube $3Q_i$ of size $3 \rho$ which contains $Q_i$ 
in the interior, satisfies
$$|3Q_i \cap E| \ge c_0|Q_i|, \quad |3Q_i \cap \CC  E| \ge c_0|Q_i|,$$
for some explicit constant $c_0>0$.
This implies that for all small $\eps$,
$$|\{u_\eps > 1-\sigma\}\cap 3Q_i| \ge \sigma |3Q_i|, \quad |\{u_\eps < 
-1+\sigma\}\cap 3Q_i| \ge \sigma |3Q_i|,$$
for some small, fixed $\sigma>0$. By Lemma \ref{help} we obtain
$$I_\eps(u_\eps,3Q_i) \ge c \rho^{n-1} \quad \mbox{if $Q_i \in F$.}$$
We write $$\bigcup_{Q_i \in F}3Q_i=\bigcup_{k=1}^N\bigcup_{Q_i \in 
F_k}3Q_i$$
with $N$ depending only on $n$ so that for each $F_k$, all cubes $3Q_i$ 
with $Q_i \in F_k$ are disjoint. We obtain
$$\sum_{Q_i \in F}I_\eps(u_\eps, 3Q_i) \le N I_\eps(u_\eps, \Omega) \le N 
C_0,$$ hence the number of cubes $Q_i$ in $F$ is bounded by $C\rho^{1-n}$.
 In conclusion $$\int_{\Omega_{\sqrt n \rho}}|\nabla \tilde \phi_\rho| \, 
dx \le C,$$ with $C$ depending on $n$, $s$ and $W$.
Since $\tilde \phi_\rho \to \chi_E-\chi_{\CC E}$ as $\rho \to 0^+$, the 
desired result follows from the lower-semicontinuity of the $BV$ norm.
\end{proof}

\section{$\Gamma$-convergence when~$s\in[1/2,1)$}\label{SS2}

In this section we prove Theorem~\ref{TH1}
and Theorem~\ref{comp} (ii) when~$s\in[1/2,1)$. In the classical case $s=1$, the $\Gamma$-convergence is obtained by relating the energy 
$\FF_\eps(u,\Omega)$ with the area of the level sets of $u$ using the coarea formula:
\begin{align*}
\int_\Omega \frac {1}{2 \eps}|\nabla u|^2+\eps W(u) \,dx &\ge \int_\Omega |\nabla u|\sqrt{2W(u)} \, dx\\ 
&=\int_{-1}^1 \sqrt{2W(s)}\mathcal H^{n-1}\big( \{u=s\} \big)\, ds.
\end{align*}    
Such formula is not available when $s<1$, so we need a careful analysis of the local and nonlocal contributions in the energy functional $\FF_\eps$.  
We will see that in the case when $s \ge 1/2$ the contribution 
$u(\Omega, \CC \Omega)$ in the kinetic term of $\FF_\eps(u, \Omega)$ 
for a minimizer $u$ becomes negligible as $\eps \to 0^+$. 

Let $D \subseteq \Omega$ be a non-empty open bounded subset of $\Omega$ 
with smooth
boundary. For all small $t >0$ define
$$D_t=\{ x \in D  : \ d_{\partial D}(x) > t \},$$
where $d_{\partial D}(x)$ represents the distance from the point $x$ to 
$\partial D$.

Next result gives an energy bound for the interpolation of two functions 
$u_k$, $w_k$ across $\partial D$: for this a fine analysis on the 
integrals is needed, which will
be accomplished by a deep modification of an argument in~\cite{DG},
where a suitable ``shell'' is selected in order to
pick up the ``right'' interpolation. Here, the situation
is much more complicate, due to the local versus
nonlocal interplay. A different nonlocal interpolation
has been recently, and independently, performed
in~\cite{Foc} in the framework of homogeneization
theory for the obstacle problem.

\begin{prop}{\label{ape3}} Fix~$\delta>0$.
Let $\varepsilon_k \to 0^+$, and let
$u_k$, $w_k$ be two sequences respectively in $L^1(D)$ and in
$L^1(\R^n)$ such that
$$u_k-w_k \to 0 \quad \mbox{as $k\rightarrow
+\infty$, in $L^1(D \setminus 
D_\delta)$.}$$
Then, there exists a
sequence $v_k$ with the following
properties:

1)
$$ v_k(x)=\left \{
\begin{array}{l}
u_k(x) \quad \quad \mbox{if $x \in D_\delta$}\\
w_k(x) \quad \quad \mbox{if $x \in \Omega\setminus D$}
\end{array}
\right.
$$

2)$$ \limsup_{k \to +\infty} \FF_{\varepsilon_k}(v_k,\Omega) \le 
\limsup_{k \to
+\infty}\Big(
\FF_{\varepsilon_k}(w_k,\Omega) -\FF_\eps(w_k, D_\delta) + 
I_{\varepsilon_k}(u_k,D) \Big) .$$
\end{prop}

\begin{proof}
Assume that there exists $C_0>0$ such that
\begin{equation}\label{11}\FF_{\varepsilon_k}(w_k,\Omega) -\FF_\eps(w_k, 
D_\delta) + I_{\varepsilon_k}(u_k,D) \le C_0,
\end{equation}
otherwise there is nothing to prove.

For simplicity of notation we drop the subindex $k$.

Since
\begin{equation}\label{BLP}\begin{split}
\KK(w, \Omega)-\KK(w, D_\delta)&=\frac{1}{2} w(\Omega \setminus 
D_\delta,\Omega \setminus D_\delta) + w(\Omega \setminus D_\delta, \CC \Omega)\\
& \ge \frac12 w(\Omega \setminus D_\delta, \CC D_\delta)
,\end{split}\end{equation}
from~\eqref{11} we obtain for $s> 1/2$ that
$$w(\Omega \setminus D_\delta, \CC 
D_\delta)+u(D \setminus D_\delta,D) \le 2C_0 \eps^{1-2s},$$
and for $s=1/2$ that
$$w(\Omega \setminus D_\delta, \CC D_\delta)+u(D 
\setminus D_\delta,D) \le 2C_0 |\log \eps|.$$
Fix $\sigma>0$ small. Let $$\tilde \delta:=\frac \delta M$$ for some large 
$M$ depending on $\sigma$, and we partition $D \setminus D_\delta$ into 
$M$ sets (i.e., ``shells")
$$D \setminus D_{\tilde \delta}, \quad D_{\tilde \delta}\setminus D_{2\tilde \delta}, \cdots,D_{(M-1)\tilde \delta}\setminus D_{M\tilde \delta}.$$
If $s > 1/2$,
\begin{align*}
2C_0 \eps^{1-2s} &\ge w(D \setminus D_\delta, \CC D_\delta)+u(D \setminus D_\delta,D)\\
 & =\sum_{j=0}^{M-1} \left( w(D_{j\tilde \delta}\setminus D_{(j+1)\tilde \delta},\CC D_\delta)+u(D_{j\tilde \delta}\setminus D_{(j+1)\tilde \delta},D) \right),
 \end{align*}
thus there exists $j \le M-1$ such that
$$w(D_{j\tilde \delta}\setminus D_{(j+1)\tilde \delta},\CC D_\delta)+u(D_{j\tilde \delta}\setminus D_{(j+1)\tilde \delta},D) \le \sigma \eps^{1-2s},
$$
provided that we choose $M$ sufficiently large. We denote
\begin{equation}\label{PRE31}
\tilde D:=D_{j\tilde \delta},\end{equation} hence, if $s> 1/2$,
we see that
\begin{equation}\label{oo1}
w(\tilde D \setminus \tilde D_{\tilde \delta}, \CC D_{\delta}) + u(\tilde D \setminus \tilde D_{\tilde \delta}, \tilde D) \le  \sigma \eps^{1-2s}.
\end{equation}
Similarly, if $s=1/2$, then
\begin{equation}\label{oo11}
w(\tilde D \setminus \tilde D_{\tilde \delta}, \CC D_{ \delta}) + u(\tilde D \setminus \tilde D_{\tilde \delta}, \tilde D) \le \sigma |\log \eps|.
\end{equation}
We remark that, since $j\le M-1$ in \eqref{PRE31}, we have that
$j\tilde\delta+\tilde\delta\le\delta$, and so
\begin{equation}\label{PRE312}
\tilde D_{\tilde\delta}\supseteq D_\delta.
\end{equation}
Next we consider $N$ shells of width $\eps\ll \tilde \delta$ of $\tilde D$,
namely
$$A_{i}:=\left \{ x\in \tilde D\,:\;i \eps < d_{\partial \tilde D} (x) \le
(i+1) \eps \right \}$$
for $0 \le i \le N-1$, with~$N$ equal the integer
part of~${\tilde \delta}/({2 
\eps})$.

We note that
\begin{equation}\label{R1}
A_i\subseteq \tilde D\setminus \tilde D_{\tilde\delta}.
\end{equation}
Also, denote by \begin{equation}\label{4.7bis}
d_i(x):=d_{\partial \tilde 
D_{i \eps}}(x).\end{equation}
Notice that
\begin{equation}\label{R2}\begin{split}
&{\mbox{for any $x\in A_i$, we have }} d_i(x)\le\eps, {\mbox{ i.e.}}\\
&\qquad 1=\min\{1,(\eps /d_i(x))^{2s}\},
\end{split}\end{equation}
while
\begin{equation}\label{R3}\begin{split}
&{\mbox{for any $x\in \tilde D_{(i+1)\eps}\setminus
\tilde D_{\tilde\delta}$, we have }} d_i(x)\ge\eps, {\mbox{ i.e.}}\\
&\qquad (\eps /d_i(x))^{2s}=\min\{1,(\eps /d_i(x))^{2s}\}.
\end{split}\end{equation}
Now, we claim that there exists $0 \le i \le N-1$ such that if $s>1/2$
\begin{equation}\label{oo2}
\int_{A_i}|u-w| \, dx + \eps^{2s}\int_{\tilde D_{(i+1) \eps}\setminus \tilde D_{\tilde \delta}}|u-w|d_i(x)^{-2s} \, dx \le
\sigma \eps,
\end{equation}
 or if $s= 1/2$
\begin{equation}\label{oo3}
\int_{A_i}|u-w| \, dx + \eps^{2s}\int_{\tilde D_{(i+1) \eps}\setminus \tilde D_{\tilde \delta}}|u-w|d_i(x)^{-2s} \, dx \le \sigma \eps |\log \eps|.
\end{equation}
Indeed, by \eqref{R1}, \eqref{R2} and \eqref{R3},
we have that the sum of all $N$ left 
hand sides for $i=0,..,N-1$ is bounded by
\begin{equation}\label{R4}2\int_{\tilde D \setminus \tilde D_{\tilde 
\delta}}|u-w| \left 
(\sum_{i=0}^{N-1}\min\{1,(\eps /d_i(x))^{2s}\}\right) \, dx.\end{equation}
Now, fix $x\in {\tilde D}\setminus{ \tilde D}_{\tilde\delta}$
and consider the shell containing $x$, that is let $i_x\in \N\cap
[0,4N]$ such that $x\in A_{i_x}$. Then, if $d_i(x)\le\eps$,
we have that $|i-i_x|\le 1$, since the other shells are
more than $\eps$ far apart from $x$. Also, if
$|i-i_x|\ge 2$, then $d_i(x)\ge (\eps/2) |i-i_x|$. {F}rom these
considerations, we see that
the sum inside the integral is bounded by a universal constant if $s > 
1/2$ or by a constant times $\log N$ if $s= 1/ 2$.

Thus the 
integral in \eqref{R4} is 
bounded by \begin{equation}\label{L1}\int_{D \setminus D_\delta}|u-w|\,dx 
, \quad \mbox {if 
$s>1/2$},\end{equation}
or \begin{equation}\label{L2}
\log N \int_{D \setminus D_\delta}|u-w|\,dx , 
\quad \mbox {if 
$s=1/2$},\end{equation}
up to multiplicative constants.

By hypothesis (for all $k$ large enough), the quantity
$$\int_{D \setminus D_\delta}|u-w|\,dx $$ can be made arbitrarily small, 
and so the claims in \eqref{oo2}
and \eqref{oo3} follow easily from \eqref{L1} and \eqref{L2}.

Now, fix a shell $A_i$ for which \eqref{oo2} or \eqref{oo3} holds. Then we 
partition $\R^n$ into five regions $P$, $Q$, $R$, $S$, $T$ where

$$P:=\tilde D_{\tilde \delta}, \qquad Q:=\tilde D_{(i+1)\eps}\setminus \tilde D_{\tilde \delta},$$
$$R:=A_{i},\qquad S:=\Omega \setminus \tilde D_{i \eps}, \qquad T=\CC \Omega.$$
Notice that
$$ Q\cup R=\tilde D_{i\eps}\setminus \tilde D_{\tilde \delta}
\subseteq \tilde D\setminus \tilde D_{\tilde\delta}$$
and, by \eqref{PRE312},
$$ R\cup S\cup T=\CC\tilde D_{(i+1)\eps}\subseteq\CC\tilde 
D_{\tilde\delta}\subseteq \CC D_\delta.$$
Therefore, \eqref{oo1} gives that
\begin{equation}\label{est w}
w(Q\cup R,R\cup S\cup T)\le\sigma \eps^{1-2s}.
\end{equation}
We choose
$$v = \phi u +(1-\phi) w $$
where $\phi$ is a smooth cutoff function with $\phi=1$ on $P \cup Q$, $\phi=0$ on $S\cup T$, and
$$\|\nabla \phi\|_{L^\infty} \le 3/ \eps.$$

Next we use \eqref{oo1}
and \eqref{oo2} and we bound $$\KK(v, \Omega)=\frac{1}{2}v(\Omega, \Omega) 
+ 
v(\Omega, \CC \Omega),$$ in terms of double integrals of $u$ and $w$. We 
consider only the case $s> 1/2$ since the only difference when $s=1/2$ is, 
as in \eqref{oo3}, the presence of an extra $|\log \eps|$ on the right hand side.

First we notice that
\begin{equation}\label{REC}
\int_{\CC B_\alpha(x)} \frac{ dy }{|x-y|^{n+2s}}
\le C \int_\alpha^{+\infty} r^{-1-2s}\,dr\le C\alpha^{-2s}
\end{equation}
for any~$\alpha>0$, and that
\begin{align}\label{o1}
v(S, S)=w(S,S),& \quad v(S,T)=w(S,T), \nonumber \\  v(P \cup Q,P\cup Q)&=u(P \cup Q,P \cup Q).
\end{align}
If $x \in P$ and $y \in R \cup S \cup T$ then
$$|x-y| \ge \tilde \delta /2 \quad \mbox{and} \quad |v(x)-v(y)|^2\le 4.$$ 
So, we use~\eqref{REC}
and we integrate the inequality
\begin{align*}
\int_{R\cup S \cup T}\frac{|v(x)-v(y)|^2}{|x-y|^{n+2s}}dy 
\le \int_{\CC B_{\tilde \delta /2}(x)}\frac{4}{|x-y|^{n+2s}}dy 
\le C \tilde \delta^{-2s},
\end{align*}
over $x \in P$ and obtain
\begin{equation}\label{o2}
v(P,R \cup S \cup T) \le C \tilde \delta^{-2s},
\end{equation}
where~$C>0$ may also depend on~$|\Omega|$.

On the other hand, recalling~\eqref{4.7bis}, we see that
if $x \in Q$ and $y \in S \cup T$ then 
\begin{equation}\label{4.17a}
|x-y|\ge d_i(x),\end{equation}  and 
\begin{equation}\label{4.17b}
|v(x)-v(y)|^2 \le 2 |u(x)-w(x)|^2+2|w(x)-w(y)|^2.\end{equation}
Thus, using~\eqref{REC} again, we deduce from~\eqref{4.17a}
that
$$\int_{S \cup T} \frac{1}{|x-y|^{n+2s}} dy 
\le C 
d_i(x)^{-2s},$$
for any~$x\in Q$, and so
we obtain, by \eqref{oo2}, \eqref{est w} and~\eqref{4.17b} that
\begin{align}\label{o3}
v(Q,S \cup T) &\le 2w(Q, S\cup T) + C \int_Q|u-w|^2d_i(x)^{-2s} \, dx \nonumber \\
&\le C \sigma \eps^{1-2s}.
\end{align}
Moreover,
if $y \in Q$ and $x \in R$ then
\begin{equation}\label{15a}
|1-\phi(x)| \le \frac 3 \eps 
d_{i+1}(x),\end{equation} and
\begin{equation}\label{15b}
|v(x)-v(y)|^2 \le 2|u(x)-u(y)|^2 + 
2|1-\phi(x)|^2 |u(x)-w(x)|^2. \end{equation}
Since, by~\eqref{REC}, we know that
$$\int_Q\frac{1}{|x-y|^{n+2s}}dy \le C 
d_{i+1}(x)^{-2s}$$ for any~$x\in R$, we 
obtain,
by \eqref{oo1}, \eqref{PRE312}, \eqref{oo2}, \eqref{15a}
and~\eqref{15b}, that
\begin{align}\label{o4}
v(Q,R) &\le 2 u(Q,R) + \frac{C}{\eps^2} \int_R |u-w|^2 d_{i+1}(x)^{2-2s} \, dx \nonumber \\
&\le 2 u(Q,R) + C \eps^{-2s} \int_R |u-w|^2 \, dx \\
& \le 2u(\tilde D \setminus \tilde D_{\tilde \delta},\tilde D \setminus \tilde D_{\tilde \delta}) + C \eps^{-2s} \int_R |u-w| \, dx \nonumber \\
& \le C \sigma \eps^{1-2s}. \nonumber
\end{align}
Similarly we find
\begin{equation}\label{o5}
v(S\cup T,R) \le C \sigma \eps^{1-2s}.
\end{equation} Furthermore,
if $x \in R $ and $y \in R$ then
\begin{align*}
|v(x)-v(y)| \le 
|w(x)-w(y)|&+|\phi(x)(u-w)(x)-\phi(y)(u-w)(y)|\\
\le |w(x)-w(y)|&+|(u-w)(x)||\phi(x)-\phi(y)|\\
& +\phi(y)|(u-w)(x)-(u-w)(y)|\\
\le 2|w(x)-w(y)|&+|u(x)-u(y)| +|(u-w)(x)||\phi(x)-\phi(y)|,
\end{align*}
hence
\begin{equation}\label{4.20a}\begin{split}
|v(x)-v(y)|^2 \le C  (|u(x)-u(y)|^2 &+|w(x)-w(y)|^2\\
&+|(u-w)(x)|^2|\phi(x)-\phi(y)|^2 ).
\end{split}\end{equation}
Also,
since $$|\phi(x)-\phi(y)| \le \min \{1,\frac{3}{\eps}|x-y|\},$$ we find
\begin{equation}\label{4.20b}\begin{split}
\int_{\R ^n}\frac{(\phi(x)-\phi(y))^2}{|x-y|^{n+2s}}dy &\le \frac{C}{\eps^2}\int_0^\eps r^{1-2s}dr + C \int_\eps^\infty r^{-1-2s} dr\\
 &\le C \eps^{-2s}.
\end{split}\end{equation}
Therefore, using \eqref{oo1}, \eqref{PRE312},
\eqref{oo2}, \eqref{est w}, \eqref{4.20a} and~\eqref{4.20b},
we can conclude 
that
\begin{align}\label{o6}
v(R,R) &\le C \left (u(R,R)+w(R,R)+\eps^{-2s}\int_R|u-w|^2dx \right )\nonumber \\
&\le C \sigma \eps^{1-2s}.
\end{align}
Also, noticing that $S\subseteq\Omega\setminus
D_\delta$,
we obtain
\begin{eqnarray*}
&&\frac12 w(S,S)+w(S,T)+\KK(w,D_\delta)\\
&=& \frac12 w(S,S)+w(S,\CC\Omega)+
\frac12 w(D_\delta,D_\delta)+\frac12 w(D_\delta,\Omega\setminus
D_\delta)\\
&&\qquad+\frac12 w(\Omega\setminus D_\delta,D_\delta)+w(D_\delta,\CC
\Omega)\\
&=& \frac12 w(S,S)
+\frac12 w(\Omega\setminus D_\delta,D_\delta)+\frac12 
w(D_\delta,\Omega)
+w(S\cup D_\delta,\CC
\Omega)\\
&\le& \frac12 w(\Omega\setminus D_\delta,S)
+\frac12 w(\Omega\setminus D_\delta,D_\delta)+\frac12
w(D_\delta,\Omega)
+w(\Omega,\CC
\Omega)\\
&\le& \frac12 w(\Omega\setminus D_\delta,\Omega)
+\frac12
w(D_\delta,\Omega)
+w(\Omega,\CC
\Omega)\\
&=&
\frac12 w(\Omega,\Omega)+w(\Omega,\CC
\Omega)\\
&=&\KK(w,\Omega).
\end{eqnarray*}
As a consequence, observing
that~$P\cup Q\subseteq D$,
and making use of
\eqref{BLP}, \eqref{o1}, \eqref{o2},
\eqref{o3}, \eqref{o4}, \eqref{o5} and~\eqref{o6}, we find
\begin{align}\label{o7}
\KK(v, \Omega) \le & \frac{1}2 w (S,S) +w(S,T) + \frac12 u(P \cup Q, P\cup 
Q) \nonumber \\
&+ C\left(\sigma \eps^{1-2s} + \tilde \delta^{-2s}\right) \\
\le& \KK (w,\Omega) - \KK(w, D_\delta) + \frac12 u(D,D) + C\left(\sigma \eps^{1-2s} + \tilde \delta^{-2s}\right). \nonumber
\end{align}
Also,
if $x \in R$ then $$W(v) \le W(w) + C|v-w| \le W(w)+C|u-w|,$$ hence 
\eqref{oo2} gives that
\begin{align}\label{o8}
\int_\Omega W(v) &\le \int_{P \cup Q}W(u) + \int_{R\cup S}W(w) + \int_R|u-w| \, dx\nonumber \\
& \le \int_DW(u) + \int_{\Omega \setminus D_\delta} W(w) + \sigma \eps.
\end{align}

{F}rom \eqref{o7}
and \eqref{o8} we obtain (for all $\eps=\eps_k$ small)
$$\FF_\eps(v, \Omega) \le \FF_\eps(w,\Omega)-\FF_\eps(w, D_\delta) + I_\eps(u,D) + C \left (\sigma +\eps^{2s-1}\tilde \delta^{-2s} \right),$$
where $C$ depends only on~$|\Omega|$, $n$ and $s$. We remark that when 
$s=1/2$, the 
last term becomes $C(\sigma +\tilde \delta^{-1} / |\log \eps|)$. Since $\sigma$ is arbitrary the proof is complete.
\end{proof}

We recall the following result about the one-dimensional minimizer which 
is proved in \cite{PSV} (see, in particular, Theorems~1.2
and~1.3 there).

\begin{thm}\label{t1}
There exists a unique (up to translations and rotations) nontrivial global 
minimizer $u_0$ of the energy $$\EE(u, \Omega):=\KK(u,\Omega) + 
\int_\Omega W(u) \, dx,$$ which depends only on one variable. If the function $
u_0$ depends only on $x_n$, then $u_0 \in C^{1,s}$ is increasing in $x_n$ and
\begin{equation}\label{EST51}
1-|u_0(x_n)| \le C |x_n|^{-2s}, \quad \quad  |u_0'(x_n)| \le C|x_n|^{-1-2s}.
 \end{equation}There exists a constant $b_\star>0$ depending only on $s,n$ and $W$ such that as $R \to \infty$

a) if $s \in(0, 1/2)$ then 
\begin{eqnarray*} &&
\liminf_{R\rightarrow+\infty} \frac1{R^{n-2s}} \int_{B_R}
\int_{\CC B_R}\frac{|u_0(x)-u_0(y)|^2}{|x-y|^{n+2s}}\,dx\,dy\,>\,0
\\ {\mbox{and }}&&\limsup_{R\rightarrow+\infty}
\frac1{R^{n-2s}} \int_{B_R}
\int_{\CC B_R}\frac{|u_0(x)-u_0(y)|^2}{|x-y|^{n+2s}}\,dx\,dy
<+\infty;\end{eqnarray*}

b) if $s=1/2$ then $$\frac{ \EE(u_0,B_R)}{R^{n-1} \log R} \to b_\star, 
\quad \mbox{and}  \quad \frac{u_0(B_R,\CC B_R)}{R^{n-1}\log R} \to 0;$$

c) if $s\in (1/2, 1)$ then $$\frac{\EE(u_0,B_R)}{R^{n-1}} \to 
b_\star,\quad \mbox{and}  \quad \frac{u_0(B_R,\CC B_R)}{R^{n-1}}\to 0.$$

\end{thm}

Theorem~\ref{t1}
says that, as $R$ gets larger and larger, the 
contribution in $\KK(
u_0,B_R)$ from $\CC B_r$ becomes negligible if $s \ge 1/2$, however when 
$s< 1/2$ this does not happen.

The energy $\FF_\eps$ is a rescaling of the energy $\EE$ in the sense that if $u$ is defined in $\R^n$ and $u_\eps(x):=u(x/ \eps),$ then
$$\FF_\eps(u_\eps,B_\rho)=\left \{
\begin{matrix}
\eps^{n-2s} \,\EE(u,B_{\rho/ \eps}) \quad \mbox{if $s<1/2$}\\
\ \\
\displaystyle\frac{\eps^{n-1}}{|\log \eps|}\, \EE(u,B_{\rho/ \eps}) \quad 
\mbox{if $s=1/2$} \\
\ \\
\eps^{n-1}\,\EE(u,B_{\rho/ \eps}) \quad \mbox{if $s > 1/2$.}
\end{matrix}
\right.
$$
Hence if
\begin{equation}\label{res}
w_\eps(x):=u_0( x / \eps),
\end{equation}
denotes the rescaling of the one-dimensional solution $u_0$, then $w_\eps$ is a global minimizer of $\FF_\eps$. Moreover, Theorem \ref{t1} can be stated in terms of $w_\eps$ and $\FF_\eps$ as
$$ \lim_{\eps \to 0^+} \FF_\eps(w_\eps,B_\rho)=b_\star \rho^{n-2s} > 
\lim_{\eps \to 0^+} I_\eps(w_\eps,B_\rho)  \quad \mbox{if $s < 1/2$},$$
and
\begin{equation}\label{Ieps}
\lim_{\eps \to 0^+} \FF_\eps(w_\eps,B_\rho)=\lim_{\eps \to 0^+} 
I_\eps(w_\eps,B_\rho)=c_\star \omega_{n-1}\rho^{n-1}, \quad \mbox{if $s \ge 1/2$,}
\end{equation}
where 
\begin{equation}\label{DEFCSTAR}c_\star:=\frac{b_\star}{\omega_{n-1}}.\end{equation}

As a consequence of Proposition \ref{ape3} and Theorem \ref{t1} 
we obtain the following

\begin{prop}\label{red.bix.2} Let~$\alpha>0$.
If~$u_\eps$ is a sequence of functions that satisfies
\begin{equation}\label{red.bix}
\lim_{\eps\rightarrow0^+}
\int_{B_\rho}
|u_\eps(x)-{\rm sign}(x_n)|\,dx
\le \alpha \rho^{n},\end{equation}
for some $\rho>0$, then
\begin{equation*}
\liminf_{\eps\rightarrow0^+}I_\eps(u_\eps, B_\rho)
\ge \omega_{n-1}\rho^{n-1}(c_\star - \eta(\alpha)).
\end{equation*}
with $\eta(\alpha)$ depending on $\alpha$ (and $n$, $s$ and $W$) and
\begin{equation}\label{alp}
\lim_{\alpha\rightarrow 0^+}\eta(\alpha)=0.\end{equation}
\end{prop}

\begin{proof}
First we prove the statement in the particular case $\rho=1$.

Assume by contradiction that the statement fails. Then we can find a sequence of functions $u_\eps$ such that
 $$\lim_{\eps\rightarrow0^+}
\int_{B_1}
|u_\eps(x)-{\rm sign}(x_n)|\,dx
= 0,$$
and  
\begin{equation}\label{C317}\limsup_{\eps\rightarrow0^+}I_\eps(u_\eps, 
B_1)
\le \omega_{n-1}c_\star - \mu,\end{equation}
for some small $\mu>0$.

Let $w_\eps$ be defined by \eqref{res}.
Then $w_\eps$ is a global minimizer for $\FF_\eps$ i.e.
\begin{equation}\label{b1}
\FF_\eps(w_\eps,B_1) \le \FF_\eps(v_\eps,B_1)
\end{equation} for any $v_\eps$ that coincides with $w_\eps$ outside $B_1$. Since $$\int_{B_1}
|w_\eps(x)-{\rm sign}(x_n)|\,dx
\rightarrow 0 \quad \mbox{as $\eps \to 0^+$,}$$  we can apply
Proposition \ref{ape3} for $u_\eps$ and $w_\eps$ with $D=\Omega=B_1$ and 
obtain
\begin{equation}\label{b2}
\limsup_{\eps\to0^+} \FF_\eps(v_\eps,B_1) \le \limsup_{\eps\to0^+} \left 
(\FF_\eps(w_\eps,B_1)-\FF_\eps(w_\eps,B_{1-\delta})+I_\eps(u_\eps,B_1)\right).
 \end{equation}
On the other hand, by \eqref{Ieps} $$\lim_{\eps\to0^+}  
\FF_\eps(w_\eps,B_1)=\omega_{n-1}c_\star, \quad  \lim_{\eps\to0^+}  
\FF_\eps(w_\eps,B_{1-\delta})=(1-\delta)^{n-1}\omega_{n-1}c_\star,$$
hence, by \eqref{b1} and \eqref{b2}
$$ (1-\delta)^{n-1}\omega_{n-1}c_\star \le \limsup_{\eps\to0^+}  
I_\eps(u_\eps,B_1), $$
and we reach a contradiction with \eqref{C317} by choosing $\delta$ 
sufficiently small.

For the general case we define $\tilde u_{\tilde \eps}$ in $B_1$ as 
$$\tilde u_{\tilde \eps}(x):=u_\eps(\rho x).$$ Then $\tilde u_{\tilde 
\eps}$ satisfies the hypothesis above in $B_1$ with $\tilde \eps:= \eps / 
\rho$ and the result follows by scaling since
\begin{equation*} 
I_\eps(u_\eps,B_\rho)=\rho^{n-1}I_{\tilde \eps}(\tilde u_{\tilde \eps}, 
B_1) \quad \quad \quad \quad \mbox{if $s> 1/2$,}$$  and 
$$I_\eps(u_\eps,B_\rho)=\rho^{n-1}\frac{|\log(\tilde \eps)|}{|\log \eps|} 
\,I_{ \tilde \eps}(\tilde u_{\tilde \eps}, B_1) \quad \quad \quad \mbox{if 
$s= 1/2$.}\qedhere\end{equation*} \end{proof}

\subsection{Reduced boundary analysis}

The idea now is to
consider any~$u_\eps$ approaching~$\chi_E-\chi_{\CC E}$,
with~$E$ of finite perimeter. Then~\eqref{red.bix} holds, suitably scaled,
near the reduced boundary of~$E$, that will be denoted,
as usual, by ~$\partial^* E$. We refer
to~\cite{Giusti}
for the basics of the theory of sets with finite perimeter
and the definition of the reduced boundary.

Precisely, let~$\nu(p)$ denote the measure theoretic unit inner
normal
at any $p\in \partial^* E$ (see Definitions 3.3
and 3.6 of \cite{Giusti}). Then, \eqref{red.bix} holds true in small balls:

\begin{corollary}\label{Us6}
Let $E$ be a set of finite perimeter, with $0\in \partial^* E$ and
\begin{equation}\label{Tmi2}\nu(0)=e_n.\end{equation}
Suppose that, as $\eps\rightarrow0^+$,
$u_\eps$ converges to $\chi_E -\chi_{\CC E}$ in~$L^1_{\rm loc} (\R^n)$.

Then,
for any $\alpha>0$ there exists $\rho(\alpha)>0$ (depending also on $n$, $s$ and $E$) such that
if $\rho\in \big(0,\,\rho(\alpha)\big]$,
we have that
$$ \lim_{\eps\rightarrow0^+} \int_{B_\rho}
|u_\eps(x)-\,{\rm sign}\,(x_n)|\,dx \le \alpha \rho^n.$$
\end{corollary}

Corollary \ref{Us6} is a consequence of the following known property of 
$\partial^* 
E$: $$\lim_{\rho \to 0^+} \rho^{-n}\int_{B_\rho}|\chi_E-\chi_{\CC E} -{\rm 
sign}(x_n)| =0.$$

\subsection{Bounding the energy from below}

We are now in the position
of obtaining a lower bound for the energy with respect to
the perimeter of the asymptotic interface for $s\in[1/2,\,1)$, and
thus proving Theorem \ref{TH1} (i).

\begin{prop}\label{7477}
Suppose that, as $\eps\rightarrow0^+$,
$u_\eps$ converges to $\chi_E -\chi_{\CC E}$ in~$L^1_{\rm loc} (\R^n)$.
Then,
$$ \liminf_{\eps\rightarrow0^+} \FF_\eps(u_\eps,\Omega)\ge
c_\star\,{\rm Per}(E,\Omega)
.$$
\end{prop}

\begin{proof}
From Proposition \ref{7477.4} (see Section \ref{4}), we may assume that $E$ has finite perimeter in $\Omega$. Then
by Theorem 4.4 of \cite{Giusti}, we have
$$ {\rm Per}(E,\Omega)={\mathcal{H}}^{n-1} (\partial^* E\cap \Omega).$$
Consequently, by fixing $\alpha>0$, we can find a collection of
balls $\{ B_j\}_{j\in\N}$ centered at
points of $\partial^* E$ and of radius $\rho_j>0$,
conveniently small in dependence of $\alpha$, such that
\begin{equation}\label{SUi1}
{\rm Per}(E,\Omega)\le \alpha+ \omega_{n-1}\sum_{j=0}^{+\infty}
\rho_j^{n-1}.\end{equation}
In fact, we can take the above balls disjoint, because of
the Vitali's Covering Theorem (see, e.g.,~\cite{Fal}), thus
\begin{equation}\label{SUi2}
\FF_\eps(u_\eps, \Omega) \ge I_\eps(u_\eps,\Omega)
\ge \sum_{j=0}^{+\infty} I_\eps(u_\eps,B_j).
\end{equation}
Also, Corollary \ref{Us6} makes \eqref{red.bix}
hold, and so we can use
Proposition \ref{red.bix.2}
in any of these balls $B_j$. Hence,
we obtain
$$ \liminf_{\eps\to0^+}
\FF_\eps(u_\eps, \Omega) \ge \omega_{n-1}(c_\star 
- 
\eta(\alpha)) \sum_{j=0}^{+\infty} \rho_j^{n-1} \ge (c_\star - 
\eta(\alpha))({\rm Per}(E,\Omega)-\alpha),$$
and the desired result follows
by letting $\alpha \to 0^+$.
\end{proof}

\subsection{Bounding the energy from above}

Now we prove part (ii) of Theorem \ref{TH1}.

\begin{prop}\label{upbd}Let $\Omega$ be a bounded domain with Lipschitz boundary. Given a set $E$, there exists a sequence $u_\eps$ converging in $L^1(\Omega)$ to $\chi_E - \chi_{\CC E}$ such that
$$ \limsup_{\eps\rightarrow0^+}\FF_\eps (u_\eps,\Omega)\le
c_\star \,{\rm Per}
\, (E,\Omega).$$
\end{prop}

\begin{proof} It was proved in \cite{MM} that there exist open sets with smooth boundaries which approximate $E$ in $\Omega$. Precisely, given any $\sigma>0$, there exists $A$ open with $\partial A$ smooth, such that
$$ \| \chi_{A \cap \Omega}- \chi_E \|_{L^1(\Omega)} \le \sigma, \quad
\quad P(A, \Omega) \le P(E, \Omega) + \sigma,$$
$$\mathcal{H}^{n-1}(\partial A \cap \partial
\Omega)=0.$$

This shows that it suffices to prove the theorem with $A$ instead of $E$. Fix $\alpha>0$ small.

Let~$d(x)$
be the signed distance of~$x$ to~$\partial A$
with the convention that~$d(x)\ge 0$ if~$x\in A$
and~$d(x)\le 0$ if~$x\in \CC A$.

We define
$$ u_\eps(x):=u_0\left( \frac{ d(x)}{\eps}\right),$$where $u_0:\R \to [-1,1]$ is the profile of the one-dimensional minimizer of $\EE$ (see Theorem \ref{t1}).

Let us take a finite overlapping
family of balls~$\{ B_{\rho_j}(x_j)\}_{j\in\N}$
centered
at~$x_j\in\partial A$, with~$\sup_{j\in\N}\rho_j\le\alpha$,
such that
$$ \partial A\cap \overline\Omega\subseteq
\bigcup_{j=0}^{+\infty} B_{\rho_j}(x_j)$$
and
$$ {\rm Per} (A,\overline\Omega)+\alpha
\ge \omega_{n-1} \sum_{j=0}^{+\infty}\rho_j^{n-1}.$$
By compactness, we may suppose that
$$ \partial A\cap \overline\Omega\subseteq V:=
\bigcup_{j=0}^{N} B_{\rho_j}(x_j),$$
for a suitable~$N\in\N$.
Notice that
$$ \delta:= \inf_{x\in\Omega
\setminus V} |d(x)|>0.$$
Recalling~\eqref{EST51}, we have that
\begin{align}\label{EST52}
\int_{\Omega
\setminus V} W(u_\eps(x))\,dx &\le C \int_{\Omega
\setminus V} |u_\eps(x)-1|\,dx \nonumber \\
&\le  C \int_{\Omega
\setminus V} \left| \frac{d(x)}{\eps}\right|^{-2s} \,dx  \le C(\delta) \eps^{2s},
\end{align}
thus, the contribution in $\FF_\eps(u_\eps)$ from the potential energy in 
$\Omega \setminus V$ tends to $0$ as $\eps \to 0^+$.
Moreover if $|d(x)| \ge \delta/2$ then we use \eqref{EST51} and obtain
$$|\nabla u_\eps(x)| = \left |
u_0'\left (\frac{d(x)}{\eps}\right ) \right |\frac 1 \eps \le C(\delta).$$

If $x \in \Omega \setminus V$ then
$$|u_\eps(x)-u_\eps(y)| \le C(\delta)|x-y|, \quad \mbox{if $|x-y|\le \delta/2$},$$ thus
$$\int_{\R^n}\frac{|u_\eps(x)-u_\eps(y)|^2}{|x-y|^{n+2s}}\,dy \le C(\delta)\left(\int_0^{\delta /2}r^{1-2s}+\int_{\delta / 2}^\infty r^{-1-2s} dr \right) \le C(\delta).$$
We find $$u_\eps(\Omega \setminus V,\R^n) \le C(\delta),$$
which together with \eqref{EST52} gives,
$$\limsup_{\eps\to0^+} \FF_\eps(u_\eps, \Omega) \le \limsup_{\eps\to0^+}  
\FF_\eps(u_\eps,V) \le \limsup_{\eps\to0^+} \sum_{j=0}^{+\infty} 
\FF_\eps(u_\eps,B_{\rho_j}).   $$

Now we estimate each term $\FF_\eps(u_\eps,B_{\rho_j})$. We will denote by~$\eta_i(\alpha)$ suitable
functions depending only on $\alpha$, $n$, $s$ and $A$ satisfying
$$ \lim_{\alpha\rightarrow0^+}\eta_i(\alpha)=0.$$
If~$\alpha$ is small enough,
then for any~$ B_{\rho_j}(x_j)$ there exists a
diffeomorphism $$x\in B_{\rho_j}(x_j) \longrightarrow z(x) \in U_j \quad \mbox{with} \quad z_n=d(x),$$
$$|D_xz-I| \le \eta_0(\alpha), \quad U_j \subset B_{1+\eta_0(\alpha)}.$$
Changing coordinates from $x$ to $z$ we find
\begin{align*}
 \FF_\eps(u_\eps,B_{\rho_j}(x_j)) &\le (1+ \eta_1(\alpha))\FF_\eps(w_\eps,U_j)\\
  &\le (1+ \eta_1(\alpha))\FF_\eps(w_\eps ,B_{1+\eta_0(\alpha)}),
\end{align*}
where $w_\eps(z)=u_0(z_n / \eps)$. {F}rom Theorem \ref{t1},
$$\limsup_{\eps\to0^+} \FF_\eps(u_\eps,B_{\rho_j}(x_j)) \le 
(1+\eta_2(\alpha))c_\star \omega_{n-1}\rho_j^{n-1},$$
and the desired result follows by letting $\alpha \to 0^+$.
\end{proof}

We denote
\begin{equation}\label{A21}
{\rm Per}(E,\overline{U}):=\lim_{\delta\rightarrow0^+}
{\rm Per}(E,U^\delta),\end{equation}
where~$$U^\delta:=\{x \in \R^n
{\mbox{ s.t. }}
\,{\rm dist}\,(x,U) \le \delta\}.$$ 
Notice that the
limit in~\eqref{A21} exists by Monotone Convergence Theorem. 

Next we prove part (ii) of Theorem \ref{comp}.

\begin{prop}\label{7477.2}
Let $\Omega$ be a bounded open subset of $\R^n$.
Suppose that~$u_\eps$ minimizes~$\FF_\eps$ in~$\Omega$
and that, as $\eps\rightarrow0^+$,
$u_\eps$ converges to $\chi_E -\chi_{\CC E}$ in~$L^1 (\Omega)$,
for some measurable $E\subseteq\Omega$.

Then $E$ has minimal perimeter in $\Omega$ and for any open set $U\subset\subset\Omega$, we have that
\begin{equation}\label{CPR}\limsup_{\eps\rightarrow0^+}\FF_\eps 
(u_\eps,{U})\le
c_\star \,{\rm Per}
\, (E,\overline{U}).\end{equation}
\end{prop}

\begin{proof}

Let $U\subset\subset\Omega$ have smooth boundary and $\delta$ be small so that $U^\delta \subset \Omega$.

Let $F$ be a measurable set in $\Omega$ such that $F$ and $E$ coincide
outside $U$. By Propositions \ref{upbd} and \ref{7477}, there exists a sequence $w_\eps\in
L^1(U^\delta)$ which converges to $\chi_F-\chi_{\CC F}$ such that
$$\lim_{\eps \to 0^+} \FF_\eps(w_\eps, U^\delta)=c_\star {\rm 
Per}(F,U^\delta).$$
From Proposition \ref{ape3} we construct a sequence $v_\eps$ which coincides with $w_\eps$ in $U$ and with $u_\eps$ in $\CC U^\delta$ such that
such that
$$\limsup_{\eps\to0^+} \FF_\eps (v_\eps, \Omega) \le
\limsup_{\eps\to0^+} \left (\FF_\eps (u_\eps, \Omega) - \FF_\eps (u_\eps, 
U) + \FF_\eps (w_\eps, U^\delta) \right).$$
Since $u_\eps$ is a minimizer,
$$\FF_\eps(u_\eps, \Omega) \le \FF_\eps(v_\eps, \Omega),$$
hence
$$ \limsup_{\eps\to0^+} \FF_\eps (u_\eps, U) \le c_\star {\rm 
Per}(F,U^\delta).$$
We let $\delta \to 0^+$ and use Proposition \ref{7477} to find
\begin{equation}\label{emin}
c_\star{\rm Per}(E,U) \le \limsup_{\eps\to0^+} \FF_\eps (u_\eps, U) \le 
c_\star{\rm Per}(F,\overline U).
\end{equation}
Since this inequalities are valid if we replace $U$ with $U^\delta$ for 
all small $\delta$, we can conclude that $${\rm Per}(E,\Omega) \le {\rm 
Per}(F,\Omega),$$ i.e. $E$ has minimal perimeter in $\Omega$. Also, by 
taking $F=E$ in \eqref{emin} we obtain \eqref{CPR}
for smooth subsets $U$. Now the general case follows easily by approximating $U$ with smooth domains from the exterior.
\end{proof}

\end{document}